\documentclass{article}
\usepackage{mathptmx} 
\usepackage{helvet}    
\usepackage{courier}   
\usepackage{type1cm}       
\usepackage{makeidx}        
\usepackage{graphicx}   
\usepackage{subcaption}
\usepackage{multicol}   
\usepackage{amsmath}
\usepackage{amsthm}
\usepackage{amssymb}
\usepackage{mathtools}
\usepackage[bottom]{footmisc}
\usepackage{verbatim}

\makeindex         
\newcommand{\M}{\mathcal{M}}

\newcommand{\dexp}{\mathrm{dexp}}

\newcommand{\g}{\mathfrak{g}}

\newcommand{\tol}{\mathrm{tol}}
\newcommand{\infgen}{\ensuremath{\psi_*}}
\newcommand{\aux}[1]{\tilde{#1}}

\begin{document}


\title{Dynamics of the N-fold Pendulum in the framework of Lie Group Integrators}
\author{Elena Celledoni \and Ergys Çokaj \and Andrea Leone \and Davide Murari \and Brynjulf Owren}

%
%
\maketitle
\abstract{Since their introduction, Lie group integrators have become a method of choice in many application areas. Various formulations of these integrators exist, and in this work we focus on Runge--Kutta--Munthe--Kaas methods. First, we briefly introduce this class of integrators, considering some of the practical aspects of their implementation, such as adaptive time stepping. We then present some mathematical background that allows us to apply them to some families of Lagrangian mechanical systems. We conclude with an application to a nontrivial mechanical system: the N-fold 3D pendulum.}
\section{Introduction}\label{introduction}
Lie group integrators are used to simulate problems whose solution evolves on a manifold. Many approaches to Lie group integrators can be found in the literature, with several applications for mechanical systems (see, e.g. \cite{celledoni2021},\cite{iserles2000}, \cite{celledoni2014}).

The present work is 
motivated by applications in modelling and simulation of slender structures like beams, and the example considered here is a chain of pendulums. The dynamics of this mechanical system is described in terms of a Lie group $G$ acting transitively on the phase space $\mathcal{M}$. This setting is used to build also a numerical integrator. 

In Section \ref{sec:liegroup} we give a brief overview of the Runge--Kutta--Munthe--Kaas (RKMK) methods with particular focus on the variable step size methods, which we use later in subsection \ref{subsec:experiments} for the numerical experiments.

In Section \ref{sec:math} we introduce some necessary mathematical background that allows us to apply RKMK methods to the system of interest. In particular, we focus on a condition that guarantees the homogeneity of the tangent bundle $TQ$ of a manifold $Q$. We then consider Cartesian products of homogeneous manifolds.

In Section \ref{sec:pendulum} we reframe the ODE system of the chain of $N$ connected $3D$ pendulums in the geometric framework presented in Section \ref{sec:math}. We write the equations of motion and represent them in terms of the infinitesimal generator of the transitive action. 
The final part shows some numerical experiments where the constant and variable step size methods are compared.

\section{RKMK methods with variable step size}\label{sec:liegroup}
The underlying idea of RKMK methods is to express a vector field $F\in\mathfrak{X}(\M)$ as $F\vert_m = \psi_*(f(m))\vert_m$, where $\psi_*$ is the infinitesimal generator of $\psi$, a transitive action on $\M$, and $f:\M\rightarrow\mathfrak{g}$. This allows us to transform the problem from the manifold $\M$ to the Lie algebra $\g$, on which we can perform a time step integration. We then map the result back to $\M$, and repeat this up to the final integration time.  More explicitly, let $h_n$ be the size of the $n-$th time step, we then update $y_n\in\M$ to $y_{n+1}$ by
\begin{equation}\label{eq:update}
    \begin{cases}
    \sigma(0) = 0\in\mathfrak{g},\\
    \dot{\sigma}(t) = \dexp_{\sigma(t)}^{-1}\circ f\circ \psi (\exp(\sigma(t)),y_n)\in T_{\sigma(t)}\mathfrak{g}, \\
    y_{n+1} = \psi(\exp(\sigma_1),y_n)\in \M,
    \end{cases}
\end{equation}
where $\sigma_1\approx \sigma(h_n)\in\mathfrak{g}$ is computed with a Runge-Kutta method.

One approach for varying the step size is based on embedded Runge--Kutta pairs for vector spaces. This approach consists of a principal method of order $p$, used to propagate the numerical solution, together with some auxiliary method, of order $\tilde{p}<p$, that is only used to obtain an estimate of the local error. This local error estimate is in turn used to derive a step size adjustment formula that attempts to keep the local error estimate approximately equal to some user-defined tolerance $\tol$ in every step.
Both methods are applied to solve the ODE for $\sigma(t)$ in (\ref{eq:update}), yielding two approximations $\sigma_1$ and $\tilde{\sigma}_1$ respectively, using the same step size $h_n$. Now, some distance measure between $\sigma_1$ and  $\tilde{\sigma}_1$ provides an estimate $e_{n+1}$ for the size of the local truncation error. Thus,
$e_{n+1}=C
h_{n+1}^{\aux{p}+1}+\mathcal{O}(h^{\aux{p}+2})$. Aiming at $e_{n+1}\approx\tol$ in every step, one may use a formula of the type
\begin{equation} \label{stepsizecontrol}
h_{n+1} = \theta\left(\frac{\tol}{e_{n+1}}\right)^{\tfrac{1}{\aux{p}+1}}\, h_n
\end{equation}
where $\theta$ is typically chosen between $0.8$ and $0.9$.
If $e_n>\tol$, the step is rejected. Hence, we can redo the step with the step size obtained by the same formula.

\section{Mathematical background}\label{sec:math}
This section introduces the mathematical background that allows us to study many mechanical systems in the framework of Lie group integrators and Lie group actions. In particular, we provide some results that we use to study the model of a chain of $N$ 3D-pendulums presented in the last section.
\subsection{The tangent bundle of some homogeneous manifolds is homogeneous}\label{subsec:tangent}
For Lagrangian mechanical systems, the phase space is usually the tangent bundle $TQ$ of some configuration manifold $Q$. In \cite{brockett1972} the authors present a setting in which the homogeneity of $Q$ implies that of $TQ$. We now briefly review and reframe it in the notation used throughout the paper.

Consider a smooth $n-$dimensional manifold $Q$, endowed with a transitive $G$-group action $\Lambda:G\times Q\rightarrow Q$. Assume that for each $q\in Q$, the map $\Lambda_q:G\rightarrow Q$ defined as $\Lambda_q(g):=\Lambda(g,q)$, is a submersion at $e\in G$. When these hypotheses hold, it can be shown that $TQ$ is a homogeneous manifold as well, and an explicit transitive action can be obtained from $\Lambda$. Let $\Lambda_*$ be the infinitesimal generator of the group action $\Lambda$, and denote with $\bar{\xi}(q):=\Lambda_*(\xi)(q)\in T_qQ$
the differential at the identity element $e\in G$ of $\Lambda_q$, evaluated at $\xi\in\mathfrak{g}$.
We then introduce $\Lambda_g:Q\rightarrow Q,$ $q\mapsto \Lambda(g,q)$ and call $T_{\bar{q}}\Lambda_g$ its tangent lift at $\bar{q}\in Q$.

Consider the manifold $\bar{G}:=G \ltimes \mathfrak{g}$, equipped with the semi-direct product Lie group structure (see, e.g., \cite{engo2003}). We can introduce a transitive group action on $TQ$ as follows:
\[
\varphi:\bar{G}\times TQ\rightarrow TQ,\;\;((g,\xi),(q,v))\mapsto \left(\Lambda(g,q),\bar{\xi}(\Lambda(g,q)) + T_{q}\Lambda_g(v)\right).
\]
By direct computation and basic properties of Lie groups (see, e.g., \cite{hall2015}), it can be seen that the action $\varphi$ is well defined. Since the action $\Lambda$ is transitive on $Q$, we notice that 
\[
\forall q'\in Q,\;\exists\,g\in G\text{ s.t. } q' = \Lambda(g,q).
\]
Then, since $\Lambda_q$ is assumed to be a submersion at $e\in G$, 
\[
\forall v'\in T_{q'}Q\;\;\exists \xi\in\mathfrak{g}\text{ s.t. }
\Lambda_*(\xi)(q') = \bar{\xi}(\Lambda(g,q))= v' - T_q\Lambda_g(v).
\]
Thus, we conclude that $\M = TQ$ is a homogeneous manifold.

In the application treated in the next section, we are interested in the case in which $Q=S^2\subset\mathbb{R}^3$, i.e. the unit sphere. In this setting, a transitive group action $\Lambda$ is given by
\[
\Lambda : SO(3)\times S^2\rightarrow S^2,\;\;(R,q)\mapsto Rq,
\]
\[
T_qS^2\ni \Lambda_*(\xi)(q) = \bar{\xi}(q) = \xi\times q,\quad 
T_q\Lambda_R(v) = Rv\in T_{Rq}S^2.
\]
Therefore, in this case we recover the restriction to $TS^2\subset\mathbb{R}^6\simeq \mathfrak{se}(3)$ of the Adjoint action of $\bar{G}=SE(3)= SO(3)\ltimes \mathbb{R}^3\simeq SO(3)\ltimes \mathfrak{so}(3)$ (see, e.g., \cite{holm2009})
\begin{equation}\label{eq:action}
\varphi((R,r),(q,v)) = (Rq, Rv + r\times Rq) =\footnote{Here $\hat{r} = \begin{bmatrix} 0 & -r_3 & r_2 \\ r_3 & 0 & -r_1 \\ -r_2 & r_1 & 0 \end{bmatrix}$, where $r=\begin{bmatrix} r_1 \\ r_2 \\ r_3\end{bmatrix}$.}
(Rq, Rv + \hat{r} Rq)
\end{equation}
which hence becomes a particular case of a more general framework.

\subsection{The Cartesian product of homogeneous manifolds is homogeneous}\label{subsec:cartesian}
Consider a family of homogeneous manifolds $\M_1,...,\M_n$. Call $(G_i,\odot_i)$ the Lie group acting transitively on the associated smooth manifold $\M_i$, and $\varphi_i$ such a transitive action. Let $\mathfrak{g}_i$ be the Lie algebra of $G_i$,  $i=1,...,n$, and \[\M = \M_1\times \M_2\times ...\times \M_n,\quad
G =G_1\times G_2 \times ... \times G_n.\]
The manifold $G$ can be naturally equipped with a Lie group structure given by the direct product. More precisely, for a pair of elements $G\ni g_i = (g_i^1,...,g_i^n),$  $i=1,2$, we can define their product $ g_1\cdot g_2 := (g_1^1\odot_1g_2^1,...,g_1^n\odot_ng_2^n)\in G.$
We can similarly define componentwise the exponential map.

This construction ensures that the manifold $\M$ is homogeneous too, and $G$ acts transitively on it. That is, let
\[
g = (g^1,...,g^n)\in G,\quad m = (m^1,...,m^n)\in \M,
\]
then
\[
\varphi: G\times \M\rightarrow \M, \quad \varphi(g,m):=(\varphi_1(g^1,m^1),...,\varphi_n(g^n,m^n)).
\]

We now restrict to the specific case $\M_i = TS^2$ for $i=1,...,n$. Since $TS^2$ is a homogeneous manifold with transitive action $\varphi$ defined as in equation (\ref{eq:action}), we can write the transitive group action
\[
\psi : (SE(3))^n\times (TS^2)^n\rightarrow (TS^2)^n,\]\[
\psi((g^1,...,g^n),(m^1,...,m^n))=(\varphi(g^1,m^1),...,\varphi(g^n,m^n)),
\]
where $g^i:=(R_i,r_i)\in SE(3)$, $m^i = (q_i,v_i)\in TS^2$.

\section{The N-fold 3D pendulum}\label{sec:pendulum}
We now apply the geometric setting from section \ref{sec:math} to the specific problem of a chain of $N$ connected 3D pendulums, whose dynamics evolves on $(TS^2)^N$.

\subsection{Equations of motion}\label{eqmoto}

\begin{figure}[htbp]
    \centering
    \includegraphics[width=.3\textwidth]{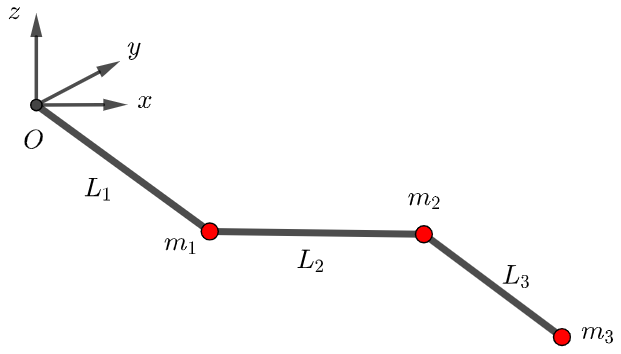}
    \caption{\footnotesize{Chain of 3 connected pendulums at a fixed time instant.}}
    \label{fig:nfold}
\end{figure}

Let us consider a  chain of $N$ pendulums subject to constant gravity $g$. The system is modeled by $N$ rigid, massless links serially connected by spherical joints, with the first link connected to a fixed point placed at the origin of the ambient space $\mathbb{R}^3$, as in figure \ref{fig:nfold}. We neglect friction and interactions among the pendulums. 

The modeling part comes from \cite{lee2018} and we omit details. We denote by $q_i\in S^2$ the configuration vector of the $i-$th mass, $m_i$, of the chain. Following \cite{lee2018}, we express the Euler–Lagrange equations for our system in terms of the configuration variables $(q_1,\dots,q_N)\in (S^2)^N\subset\mathbb{R}^{3N}$, and their angular velocities $(\omega_1,...,\omega_N)\in T_{q_1}S^2\times ... \times T_{q_N}S^2\subset\mathbb{R}^{3N}$, defined be the following kinematic equations:
\begin{equation} \label{kineqs}
    \dot{q}_i = \omega_i\times q_i, \quad i=1,\dots,N.
\end{equation}
The Euler–Lagrange equations of the system can be written as
\begin{equation} \label{eleqs}
    R(q)\dot{\omega} = \left[\sum_{\substack{j=1\\ j\neq i}}^N M_{ij}|\omega_j|^2\hat{q}_i q_j - \Big(\sum_{j=i}^N m_j\Big)gL_i \hat{q}_i e_3 \right]_{i=1,...,N} = \begin{bmatrix}r_1\\ \vdots \\ r_N \end{bmatrix}\in\mathbb{R}^{3N},
\end{equation}
where $R(q)\in\mathbb{R}^{3N\times 3N}$ is a symmetric block matrix defined as

\[
R(q)_{ii} = \Big(\sum_{j=i}^Nm_j\Big)L_i^2I_3\in\mathbb{R}^{3\times 3},
\]
\[
R(q)_{ij} = \Big(\sum_{k=j}^N m_k\Big)L_iL_j\hat{q}_i^T\hat{q}_j\in\mathbb{R}^{3\times 3} = R(q)_{ji}^T,\; i<j,
\]

and 
\[
M_{ij} =\Big(\sum_{k={\text{max}}\{i,j\}}^N m_k\Big)L_iL_j I_3\in\mathbb{R}^{3\times 3}.
\]
Equations (\ref{kineqs})-(\ref{eleqs}) define the dynamics of the N-fold pendulum, and hence a vector field $F\in\mathfrak{X}((TS^2)^N)$. We now find a function $f:(TS^2)^N\rightarrow \mathfrak{se}(3)^N$ such that
\[
\infgen(f(m))\vert_m = F\vert_m,\;\;\forall m\in (TS^2)^N,
\]
where $\psi$ is defined as in subsection \ref{subsec:cartesian}.

Since $R(q)$ defines a linear invertible map (see \cite{celledoni2021})
\[
A_{q}:T_{q_1}S^2\times ... \times T_{q_N}S^2 \rightarrow T_{q_1}S^2 \times ... \times T_{q_N}S^2,\quad A_q(\omega):=R(q)\omega,
\]
we can rewrite the ODEs for the angular velocities as follows:
\begin{equation}\label{eq:omegadot}
\dot{\omega}= A_{q}^{-1}\left(\begin{bmatrix}r_1\\ \vdots \\ r_N \end{bmatrix}\right) =\begin{bmatrix}
h_1(q,\omega) \\ \vdots \\ h_N(q,\omega)
\end{bmatrix} = \begin{bmatrix}
a_1(q,\omega)\times q_1 \\
\vdots \\
a_N(q,\omega)\times q_N
\end{bmatrix}.
\end{equation}
In equation (\ref{eq:omegadot}) the $r_i$s are defined as in (\ref{eleqs}),
and $a_1,...,a_N:(TS^2)^N\rightarrow \mathbb{R}^3$ can be defined as $a_i(q,\omega):=q_i\times h_i(q,\omega)$. Thus, the map $f$ is given by

\[
f(q,\omega) = \begin{bmatrix}
\omega_1 \\
q_1\times h_1(q,\omega) \\ \vdots \\ \omega_N \\ q_N\times h_N(q,\omega)
\end{bmatrix}\in\mathfrak{se}(3)^N\simeq \mathbb{R}^{6N}.
\]

\subsection{Numerical experiments}\label{subsec:experiments}
In this section we show a numerical experiment with the N-fold 3D pendulum, in which we compare the performance of constant and variable step size methods. We do not show results on the preservation of the geometry (up to machine accuracy), since this is given by construction. We consider the RKMK pair coming from Dormand–Prince method (DOPRI 5(4) \cite{dormand1980family}, which we denote by RKMK(5,4)). We set a tolerance of $10^{-6}$ and solve the system with the RKMK(5,4) scheme. Fixing the number of time steps required by RKMK(5,4), we repeat the experiment with RKMK of order 5 (denoted by RKMK5). The comparison occurs at the final time $T=3$ using the Euclidean norm of the ambient space $\mathbb{R}^{6N}$. The quality of the approximation is measured against a reference solution obtained with ODE45 from MATLAB with a strict tolerance.

\begin{figure}[h!]
\centering
\begin{subfigure}{.45\textwidth}
    \centering
    \includegraphics[width=\textwidth]{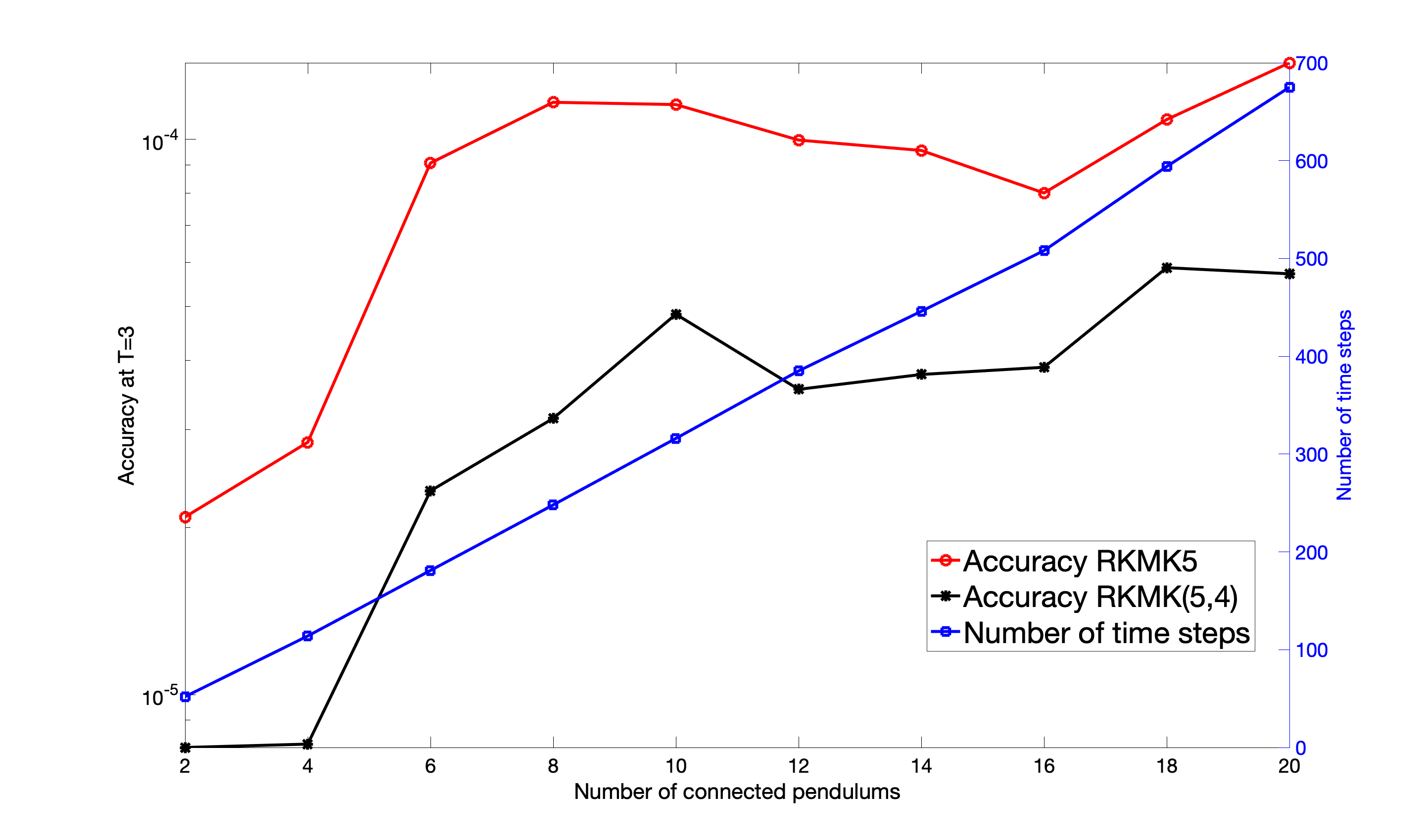}
    \caption{\scriptsize{Accuracy against the number of pendulums}}
    \label{fig:figA}
\end{subfigure}
\begin{subfigure}{.45\textwidth}
    \centering
    \includegraphics[width=\textwidth]{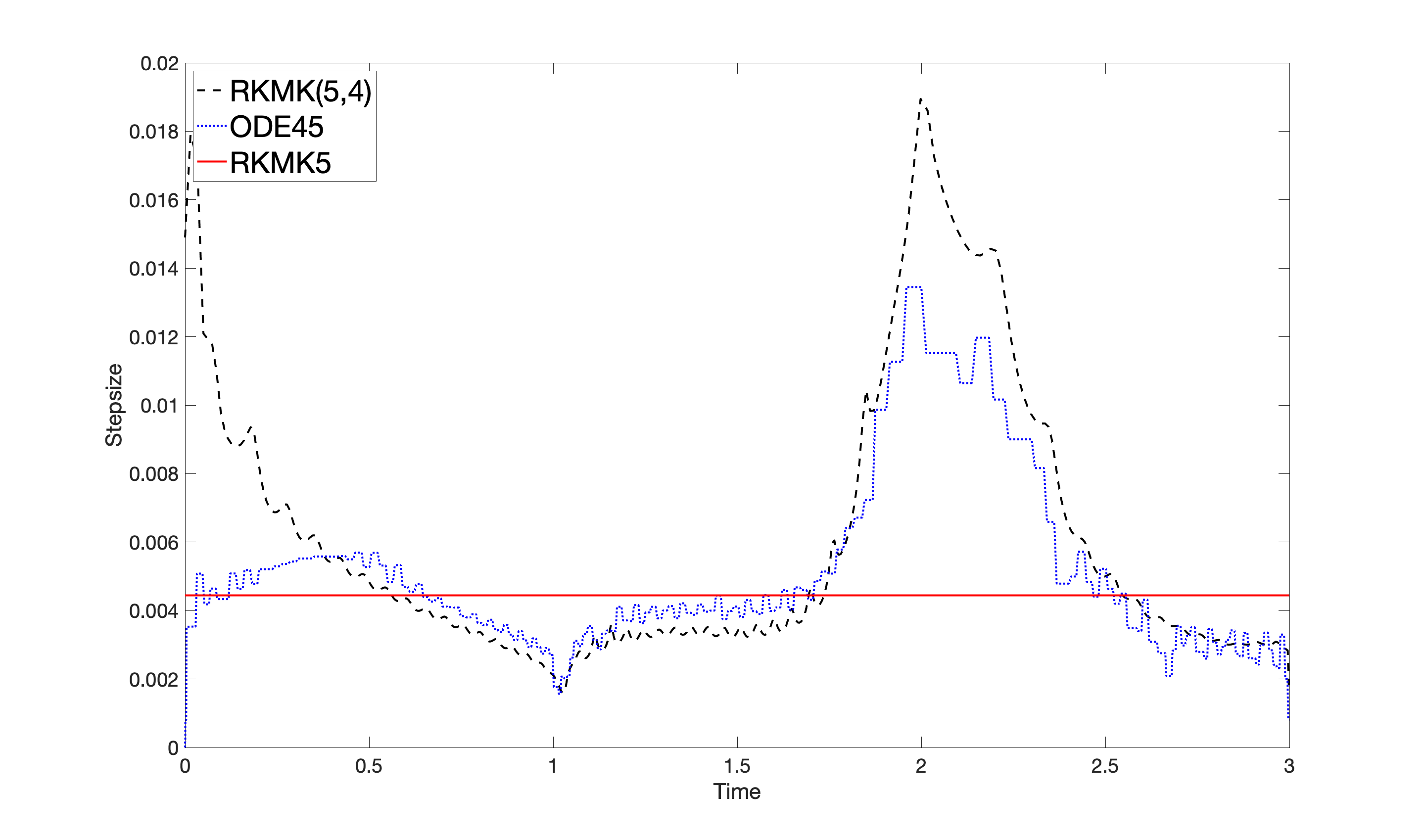}
    \caption{\scriptsize{Comparison of step sizes with $20$ pendulums}}
    \label{fig:figB}
\end{subfigure}
\caption{\footnotesize{Analysis of variable step size against constant one}}
\label{fig:variableStep}
\end{figure}

The motivating application behind the choice of this mechanical system has been some intuitive relation with flexible slender structures like beams. For this limiting behaviour to make sense, we first fix the length of the entire chain of pendulums to some $L$, then we set the size of each pendulum to $L_i=L/N$ and initialize $ (q_i,\omega_i) = \left(1,0,0,0,0,0\right)$, $\forall i=1,...,N$.  As we can see in figure \ref{fig:figA}, the results of our experiments show that number of time steps that RKMK(5,4) requires to reach the desired accuracy increases with $N$, and this can be read in terms of an augmentation of the dynamics' complexity. For this reason, as highlighted in figure \ref{fig:variableStep}, distributing these time steps uniformly in the time interval $[0, T]$ becomes an inefficient approach, and hence a variable step size method gives better performance.

\end{document}